\documentclass[12pt,a4paper]{article} 
\usepackage{amssymb}
\usepackage{times} 
\newcommand{\nc}{\newcommand} 
\nc{\bb}{\bigskip}
\nc{\C}{\mathbb{C}}
\nc{\cl}{\centerline} 
\nc{\ds}{\displaystyle} 
\nc{\ind}{\hskip 1em\relax}
\nc{\N}{\mathbb{N}}
\nc{\prodj}{\prod_{\stackrel{\scr j=1,\dots,k}{j\not=p}}}
\nc{\R}{\mathbb{R}}
\nc{\scr}{\scriptstyle} 
\nc{\sumb}{\sum_{\stackrel{\scr\beta_p=0}{|\beta|=n}}}

\begin{document}

\cl{\Large Matrix exponentials}\bb\bb

Let $A$ be a complex square matrix, put 

$$Q(z)=(z-a_1)^{\alpha_1+1}\ \cdots\ (z-a_k)^{\alpha_k+1},$$

with the $a_p$ distinct and the $\alpha_p$ nonnegative integers,  
assume $Q(A)=0$, set $$Q_p(z):=\frac{Q(z)}{(z-a_p)^{\alpha_p+1}}\quad,$$  
let $b_{p,n}$ be the $n$-th Taylor coefficient of $1/Q_p(z)$ at 
$z=a_p$, let $f$ be an entire function, and let $P(f(z),z)\in\C[z]$ be 
$Q(z)$ times the singular part of $f(z)/Q(z)$. \bb

{\bf Theorem.} {\it We have} \begin{itemize}
 
\item $\ds b_{p,n}=(-1)^n \sumb\ \prodj\ {\alpha_j+\beta_j\choose\alpha_j}
\quad\frac{1}{(a_p-a_j)^{\alpha_j+1+\beta_j}}$

{\it where $\beta$ runs over $\N^k$ and} 
$|\beta|:=\beta_1+\cdots+\beta_k$, 

\item $\ds P(f(z),z)=\sum_{p=1}^k\ \sum_{q=0}^{\alpha_p}\ \sum_{j=0}^q\quad 
\frac{f^{(j)}(a_p)}{j!}\ b_{p,q-j}\ (z-a_p)^q\ Q_p(z),$

\item $\ds f(A)=P(f(z),A).$\end{itemize}\bb 

\ind Let $Q$ be as above and $d$ its degree, and define $g_j:\R\to\C$ for 
$0\le j<d$ by 
$$P(e^{tz},z)=g_{d-1}(t)\ z^{d-1}+\cdots+g_1(t)\ z+g_0(t).$$

\ind If $u:\R\to\C$ is smooth satisfying $Q(\frac{d}{dt})\ u=h$ where 
$h:\R\to\C$ is continuous, then 
$$u(t)=\sum_{j=0}^{d-1}\ u^{(j)}(0)\ g_j(t)+\int_0^tg_{d-1}(t-y)\ h(y)\
dy.$$

\ind Assume that each $a_p$ is an eigenvalue, let $A=S+N$ 
($S$ semisimple, $N$ nilpotent) be the Jordan decomposition of $A$, 
and $A=\sum A_p=\sum(a_pE_p+N_p)$ be its spectral decomposition. Then 

$$E_p=\sum_{q=0}^{\alpha_p}\ b_{p,q}\ (A-a_p)^q\ Q_p(A),\quad%
N_p=\sum_{q=0}^{\alpha_p-1}\ b_{p,q}\ (A-a_p)^{q+1}\ Q_p(A).$$

\scriptsize Pierre-Yves Gaillard, Institut \'Elie Cartan, 
Universit\'e Nancy I, BP 239, 54506 Vand\oe uvre, France, 
gaillard@iecn.u-nancy.fr

\end{document}